\documentclass[12pt,reqno]{article}
\usepackage[usenames]{color}
\usepackage[colorlinks=true,linkcolor=webgreen,filecolor=webbrown,citecolor=webgreen]{hyperref}
\usepackage{amssymb,amsmath,amsfonts,amsthm}
\usepackage{enumerate}

\definecolor{webgreen}{rgb}{0,.5,0}
\definecolor{webbrown}{rgb}{.6,0,0}

\allowdisplaybreaks
\input{tcilatex}
\begin{document}

\begin{center}
\vskip1cm{\LARGE \textbf{New proofs of Melzak's identity}} \\[1.5cm]
{\large Ulrich Abel}\\[3mm]
\textit{Fachbereich MND}\\[0pt]
\textit{Technische Hochschule Mittelhessen}\\[0pt]
\textit{Wilhelm-Leuschner-Stra\ss e 13, 61169 Friedberg, }\\[0pt]
\textit{Germany}\\[0pt]
\href{mailto:Ulrich.Abel@mnd.thm.de}{\texttt{Ulrich.Abel@mnd.thm.de}} \\[1cm]
{\large Henry W. Gould}\\[3mm]
\textit{Department of Mathematics}\\[0pt]
\textit{West Virginia University}\\[0pt]
\textit{Morgantown, WV 26506}\\[0pt]
\textit{USA}\\[0pt]
\href{mailto:henrygou@gmail.com}{\texttt{henrygou@gmail.com}} \\[1cm]
{\large Jocelyn Quaintance}\\[3mm]
\textit{Department of Computer Science}\\[0pt]
\textit{University of Pennsylvania}\\[0pt]
\textit{Philadelphia, PA 19104}\\[0pt]
\textit{USA}\\[0pt]
\href{mailto:quaintan@temple.edu}{\texttt{quaintan@temple.edu}} \\[1cm]
\end{center}

\vspace{2cm} \vfill

{\large \textbf{Abstract.}}

\bigskip

In their recent book \cite{Quaintance-Gould-2016} on combinatorial
identities, Quaintance and Gould devoted one chapter \cite[Chapt. 7]%
{Quaintance-Gould-2016} to Melzak's identity. We give new proofs for this
identity and its generalization.

\bigskip



\bigskip

\newpage

\section{Introduction}

Melzak's identity \cite{Melzak-AMM-Problem-4458} (see also \cite[Eq. (7.1)]%
{Quaintance-Gould-2016}) states that 
\begin{equation}
f\left( x+y\right) =x\binom{x+n}{n}\sum_{k=0}^{n}\left( -1\right) ^{k}\binom{%
n}{k}\frac{f\left( y-k\right) }{x+k}\text{ \qquad }\left( x,y\in \mathbb{R},%
\text{ }n=0,1,2,\ldots \right) ,  \label{Melzak-identity}
\end{equation}%
where $x\neq -k$ $\left( k=0,1,\ldots ,n\right) $, for all algebraic
polynomials $f$ up to degree $n$. It has several interesting applications
(see \cite[Chapt. 7]{Quaintance-Gould-2016}).

In his recent book on combinatorial identities, J. Quaintance and H. W.
Gould gave an elementary proof \cite[p. 79--82]{Quaintance-Gould-2016} using
partial fraction decomposition. Furthermore, they generalized Melzak's
identity by replacing $x+k$ in the denominator $x+k$ of Eq. $\left( \ref%
{Melzak-identity}\right) $ with a product $\left( x_{0}+k\right) \left(
x_{1}+k\right) \cdots \left( x_{j}+k\right) $ of finitely many pairwise
different linear factors: 
\begin{equation}
\sum_{k=0}^{n}\left( -1\right) ^{k}\binom{n}{k}\frac{f\left( y-k\right) }{%
\prod\nolimits_{i=0}^{j}\left( x_{i}+k\right) }=\sum_{i=0}^{j}\frac{1}{x_{i}%
\binom{x_{i}+n}{n}\prod_{\nu =0,\nu \neq i}^{j}\left( x_{\nu }-x_{i}\right) }%
f\left( y+x_{i}\right) ,  \label{Melzak-generalization}
\end{equation}%
for $n=0,1,2,\ldots $ and $y\in \mathbb{R}$ (see \cite[Eq. (7.52)]%
{Quaintance-Gould-2016}). Of course, we assume that $x_{i}\neq -k$ for $%
k=0,1,\ldots ,n$ and $i=0,1,\ldots ,j$.

\section{An elementary proof of Melzak's identity}

Because Melzak's identity $\left( \ref{Melzak-identity}\right) $ is linear
it is sufficient to prove it for the functions $f_{r}\left( x\right) =x^{r}$%
\ $\left( r=0,1,\ldots ,n\right) $. Observing that 
\begin{equation*}
\frac{f_{r}\left( y-k\right) }{x+k}=\frac{\left( \left( x+y\right) -\left(
x+k\right) \right) ^{r}}{x+k}=\sum_{j=0}^{r}\left( -1\right) ^{j}\binom{r}{j}%
\left( x+y\right) ^{r-j}\left( x+k\right) ^{j-1}
\end{equation*}%
we obtain 
\begin{equation*}
\sum_{k=0}^{n}\left( -1\right) ^{k}\binom{n}{k}\frac{f\left( y-k\right) }{x+k%
}=\sum_{j=0}^{r}\left( -1\right) ^{j}\binom{r}{j}\left( x+y\right)
^{r-j}\sum_{k=0}^{n}\left( -1\right) ^{k}\binom{n}{k}\left( x+k\right)
^{j-1}.
\end{equation*}%
The second sum is a difference of order $n$ which vanishes when applied to
polynomials of degree $r-1\leq n-1$. Therefore, we have, for $x>0$, 
\begin{eqnarray*}
&&\sum_{k=0}^{n}\left( -1\right) ^{k}\binom{n}{k}\frac{f_{r}\left(
y-k\right) }{x+k} \\
&=&\left( x+y\right) ^{r}\sum_{k=0}^{n}\left( -1\right) ^{k}\binom{n}{k}%
\left( x+k\right) ^{-1}=f_{r}\left( x+y\right) \sum_{k=0}^{n}\left(
-1\right) ^{k}\binom{n}{k}\int_{0}^{1}t^{x+k-1}dt \\
&=&f_{r}\left( x+y\right) \int_{0}^{1}t^{x-1}\left( 1-t\right)
^{n}dt=f_{r}\left( x+y\right) \cdot B\left( x,n+1\right) .
\end{eqnarray*}%
Application of the well-known formula $B\left( x,n+1\right) =x\binom{x+n}{n}$
for the Beta function shows that identity $\left( \ref{Melzak-identity}%
\right) $ is valid, for the function $f_{r}$. Because both sides are
rational functions the restriction $x>0$ can be omitted.

\section{The proof by Parker}

Melzak presented his identity $\left( \ref{Melzak-identity}\right) $ in the
problem section of the American Mathematical Monthly in 1951. We sketch the
published solution by Parker \cite{Melzak-AMM-Solution-4458}. His proof (see 
\cite[p. 98]{Quaintance-Gould-2016}) used the Lagrange interpolation
polynomial. Let 
\begin{equation*}
\omega _{n}\left( x\right) =\prod\nolimits_{\nu =0}^{n}\left( x-x_{\nu
}\right) .
\end{equation*}%
The Lagrange interpolation polynomial of degree $n$ with respect to the
knots $x_{0},x_{1},\ldots ,x_{n}$ can be written in the form 
\begin{equation}
L_{n}\left( x\right) =\omega _{n}\left( x\right) \sum_{k=0}^{n}\frac{f\left(
x_{k}\right) }{\left( x-x_{k}\right) \omega _{n}^{\prime }\left(
x_{k}\right) }.  \label{Lagrange-Interpolation}
\end{equation}%
Obviously, we have $\omega _{n}^{\prime }\left( x_{k}\right)
=\prod\nolimits_{\nu =0,\nu \neq k}^{n}\left( x_{k}-x_{\nu }\right) $. We
choose $x_{k}=-k$. Then $\omega _{n}\left( x\right) =x\left( x+1\right)
\cdots \left( x+n\right) $ and $\omega _{n}^{\prime }\left( -k\right)
=\prod\nolimits_{\nu =0,\nu \neq k}^{n}\left( -k+\nu \right) =\left(
-1\right) ^{k}k!\left( n-k\right) !$ which leads to 
\begin{equation*}
L_{n}\left( x\right) =\binom{x+n}{n}\sum_{k=0}^{n}\left( -1\right) ^{k}%
\binom{n}{k}\frac{f\left( -k\right) }{x+k}.
\end{equation*}%
Let $f$ be a polynomial of degree $\leq n$. Then $L_{n}$ coincides with $f$
and we obtain Melzak's identity $\left( \ref{Melzak-identity}\right) $, by a
shift of the variable replacing $f\left( .\right) $ with $f\left( .+y\right) 
$.

\section{The generalization of Melzak's identity}

Now we turn to formula $\left( \ref{Melzak-generalization}\right) $. Let $%
x_{0},x_{1},\cdots ,x_{j}$ be pairwise different reals. The Lagrange
interpolation polynomial $\left( \ref{Lagrange-Interpolation}\right) $ of $%
f\left( x\right) =1$ with respect to the knots $x_{0},\ldots ,x_{j}$ is 
\begin{equation*}
1=f\left( x\right) =L_{j}\left( x\right) =\omega _{j}\left( x\right)
\sum_{i=0}^{j}\frac{1}{\left( x-x_{i}\right) \prod\nolimits_{\nu =0,\nu \neq
i}^{j}\left( x_{i}-x_{\nu }\right) }.
\end{equation*}%
Evaluation at $x=-k$ yields the partial fraction decomposition 
\begin{equation*}
\frac{1}{\prod\nolimits_{\nu =0}^{j}\left( k+x_{\nu }\right) }=\sum_{i=0}^{j}%
\frac{1}{\left( k+x_{i}\right) \prod\nolimits_{\nu =0,\nu \neq i}^{j}\left(
x_{\nu }-x_{i}\right) }
\end{equation*}%
(see \cite[Eq. (7.51)]{Quaintance-Gould-2016}). Following \cite[p. 92--93]%
{Quaintance-Gould-2016} J. Quaintance and H. W. Gould multiply the latter
equation by $\left( -1\right) ^{k}\binom{n}{k}f\left( y-k\right) $, sum over 
$k$, interchange the order of summation on the right side and apply Melzak's
identity $\left( \ref{Melzak-identity}\right) $ in order to obtain 
\begin{equation*}
\sum_{k=0}^{n}\left( -1\right) ^{k}\binom{n}{k}\frac{f\left( y-k\right) }{%
\prod\nolimits_{\nu =0}^{j}\left( k+x_{\nu }\right) }=\sum_{i=0}^{j}\frac{%
f\left( x_{i}+y\right) }{x_{i}\binom{x_{i}+n}{n}\prod\nolimits_{\nu =0,\nu
\neq i}^{j}\left( x_{\nu }-x_{i}\right) },
\end{equation*}%
for all polynomials $f$ of degree $\leq n$ (see \cite[Eq. (7.52)]%
{Quaintance-Gould-2016}). This is identity $\left( \ref%
{Melzak-generalization}\right) $.

Now we are going to prove this identity to hold even for all polynomials $f$
of degree $\leq n+j$ as it was claimed in \cite[Eq. (7.52)]%
{Quaintance-Gould-2016}. To this end we choose $f\left( x\right) =\binom{-x}{%
j}$. The Lagrange interpolation polynomial $\left( \ref%
{Lagrange-Interpolation}\right) $ with respect to the knots $x_{0},\ldots
,x_{j}$ is 
\begin{equation*}
\binom{-x}{j}=f\left( x\right) =L_{j}\left( x\right) =\omega _{j}\left(
x\right) \sum_{i=0}^{j}\frac{1}{\left( x-x_{i}\right) \prod\nolimits_{\nu
=0,\nu \neq i}^{j}\left( x_{i}-x_{\nu }\right) }\binom{-x_{i}}{j}.
\end{equation*}%
Evaluation at $x=-k$ yields 
\begin{equation*}
\sum_{i=0}^{j}\frac{1}{\left( -k-x_{i}\right) \prod\nolimits_{\nu =0,\nu
\neq i}^{j}\left( x_{i}-x_{\nu }\right) }\binom{-x_{i}}{j}=\frac{1}{\omega
_{j}\left( -k\right) }\binom{k}{j}.
\end{equation*}%
Using $\binom{-x_{i}}{j}=\left( -1\right) ^{j}\binom{x_{i}+j-1}{j}$ and $%
\omega _{j}\left( -k\right) =\left( -1\right) ^{j+1}\prod\nolimits_{\nu
=0}^{j}\left( k+x_{\nu }\right) $ we conclude that 
\begin{eqnarray*}
&&\sum_{k=0}^{n+j}\left( -1\right) ^{k}\binom{n+j}{k}\frac{f\left(
y-k\right) }{\prod\nolimits_{\nu =0}^{j}\left( k+x_{\nu }\right) }\binom{k}{j%
} \\
&=&\sum_{i=0}^{j}\frac{\left( -1\right) ^{j}}{\prod\nolimits_{\nu =0,\nu
\neq i}^{j}\left( x_{i}-x_{\nu }\right) }\binom{x_{i}+j-1}{j}%
\sum_{k=0}^{n+j}\left( -1\right) ^{k}\binom{n+j}{k}\frac{f\left( y-k\right) 
}{k+x_{i}}.
\end{eqnarray*}%
Using $\binom{n+j}{k}\binom{k}{j}=\binom{n+j}{j}\binom{n}{k-j}$ and
application of Melzak's identity $\left( \ref{Melzak-identity}\right) $
leads to 
\begin{equation*}
\binom{n+j}{j}\sum_{k=j}^{n+j}\left( -1\right) ^{k}\binom{n}{k-j}\frac{%
f\left( y-k\right) }{\prod\nolimits_{\nu =0}^{j}\left( k+x_{\nu }\right) }%
=\sum_{i=0}^{j}\frac{\left( -1\right) ^{j}}{\prod\nolimits_{\nu =0,\nu \neq
i}^{j}\left( x_{i}-x_{\nu }\right) }\binom{x_{i}+j-1}{j}\frac{f\left(
x+y\right) }{x_{i}\binom{x_{i}+n+j}{n+j}},
\end{equation*}%
for all polynomials $f$ of degree $\leq n+j$. Using $\binom{x_{i}+n+j}{n+j}=%
\binom{x_{i}+n+j}{n}\binom{x_{i}+j}{j}\binom{n+j}{n}^{-1}$ we obtain 
\begin{equation*}
\sum_{k=0}^{n}\left( -1\right) ^{k+j}\binom{n}{k}\frac{f\left( y-k-j\right) 
}{\prod\nolimits_{\nu =0}^{j}\left( k+j+x_{\nu }\right) }=\sum_{i=0}^{j}%
\frac{1}{\prod\nolimits_{\nu =0,\nu \neq i}^{j}\left( x_{\nu }-x_{i}\right) }%
\binom{x_{i}+j-1}{j}\frac{f\left( x_{i}+y\right) }{x_{i}\binom{x_{i}+n+j}{n}%
\binom{x_{i}+j}{j}}.
\end{equation*}%
The shifts $y\rightarrow y+j$ and $x_{\nu }\rightarrow x_{\nu }-j$ $\left(
\nu =0,\ldots ,j\right) $ yield 
\begin{equation*}
\sum_{k=0}^{n}\left( -1\right) ^{k+j}\binom{n}{k}\frac{f\left( y-k\right) }{%
\prod\nolimits_{\nu =0}^{j}\left( k+x_{\nu }\right) }=\sum_{i=0}^{j}\frac{1}{%
\prod\nolimits_{\nu =0,\nu \neq i}^{j}\left( x_{\nu }-x_{i}\right) }\binom{%
x_{i}-1}{j}\frac{f\left( x_{i}+y\right) }{\left( x_{i}-j\right) \binom{%
x_{i}+n}{n}\binom{x_{i}}{j}}.
\end{equation*}%
Now formula $\left( \ref{Melzak-generalization}\right) $ follows, for all
polynomials $f$ of degree $\leq n+j$, by the observation that $\left(
x_{i}-j\right) \binom{x_{i}}{j}=x_{i}\binom{x_{i}-1}{j}$.


\begin{thebibliography}{9}
\bibitem{Melzak-AMM-Solution-4458} Z. A. Melzak, V. D. Gokhale, and W. V.
Parker, \newblock Advanced Problems and Solutions: Solutions: 4458, %
\newblock {\em Amer. Math. Monthly} {\bf 60} (1953), 53--54.

\bibitem{Melzak-AMM-Problem-4458} Z. A. Melzak, \newblock Problem 4458, %
\newblock {\em Amer. Math. Monthly} {\bf 58} (1951), 636.

\bibitem{Quaintance-Gould-2016} Jocelyn Quaintance and Henry W. Gould, %
\newblock Combinatorial Identities for Stirling Numbers: The Unpublished
Notes of H. W. Gould, \newblock World Scientific Publishing, Singapore 2016.
\end{thebibliography}
\end{document}